\documentclass[numreferences]{kluwer} 
\usepackage{oldgerm}
\usepackage{mathrsfs}
\usepackage{amsfonts}
\usepackage{mathrsfs}
\usepackage{amsthm}
\input{diagrams}

\def\b{\mathbb }

\def\phi{\varphi }

\def\reel{\b R}
\def\nat{\b N}

\def\comp{\b C}
\def\ganz{\b Z}

\def\bar{\overline }
\theoremstyle{plain}
\swapnumbers
\newtheorem{theorem}{Theorem}[section]
\newtheorem{corollary}[theorem]{Corollary}
\newtheorem{lemma}[theorem]{Lemma}
\newtheorem{proposition}[theorem]{Proposition}
\theoremstyle{definition}

\newtheorem{remark}[theorem]{Remark}

\newtheorem{examples}[theorem]{Examples}
\newtheorem{example}[theorem]{Example}

\newtheorem{hypergroups}[theorem]{Commutative hypergroups}

\newtheorem{method}[theorem]{Description of the method}
\newtheorem{moments}[theorem]{Moments}
\newtheorem{series}[theorem]{The further classical series of simple groups}

\begin{document}

\begin{article}
\begin{opening}
         
\title{$SU(d)$-biinvariant random walks on  $SL(d,\b C)$ and
 their Euclidean counterparts  } 
\author{Margit \surname{R\"osler}\email{mroesler@science.uva.nl}\thanks{Partially supported by the Netherlands Organisation for Scientific Research (NWO), project nr. B 61-544}} 
\institute{Korteweg-de Vries Institute for Mathematics, University of Amsterdam, Plantage Muidergracht 24, 1018 TV Amsterdam, Netherlands}
\author{Michael \surname{Voit}\email{michael.voit@mathematik.uni-dortmund.de}}
 \institute{Fachbereich Mathematik,  Universit\"at Dortmund,
D-44221 Dortmund, Germany} 
\runningauthor{R\"osler and Voit}
\runningtitle{$SU(d)$-biinvariant random walks on $SL(d\b, C)$}

\begin{abstract}
We establish a deformation isomorphism between  the algebras of
$SU(d)$-biinvariant compactly supported
  measures on $SL(d,\comp)$ and $SU(d)$-conjugation invariant measures on 
the Euclidean space $H_d^0$ 
 of all Hermitian   $d\times d$-matrices with trace $0$.  This 
isomorphism concisely explains a close connection between
the  spectral problem for sums of Hermititan matrices  on one hand and 
the singular spectral problem for products of matrices from $SL(d,\comp)$
on the other, which has recently been observed by Klyachko \cite{Kl2}.
From this deformation we further obtain 
an explicit, probability preserving and isometric isomorphism between the 
Banach algebra of  bounded   $SU(d)$-biinvariant measures on $SL(d,\comp)$
 and a certain (non-invariant) subalgebra of  the  bounded signed measures on $H_d^0$. We demonstrate how this probability preserving isomorphism
leads to limit
theorems for  the  singular spectrum of  $SU(d)$-biinvariant random walks on
 $SL(d,\comp)$ in a simple way.
Our construction  relies on deformations of 
 hypergroup convolutions  and will be carried out 
in the general setting of  complex  semisimple Lie groups.
\end{abstract}

%\keywords{sample, \LaTeX}

\classification{Mathematics Subject Classification}
{Primary:  43A10; Secondary: 43A85, 60B15,  43A62, 22E46.}

\end{opening}

\section{Introduction}

The spectral problem for possible sums of two random Hermitian matrices 
with given spectra had been a long-standing problem, formulated as Horn's conjecture, until it was completely solved only recently by
 Klyachko, Knutson and Tao (\cite{Kl1}, \cite{KT}).
 In \cite{Kl2}, Klyachko then observed a close  connection between
$SU(d)$-biinvariant  random walks  on  $SL(d,\comp)$
and
random walks on the  additive group $H_d^0$ of
all Hermitian $d\times d$-matrices with trace 0
whose transition probabilities are
conjugation-invariant under  $SU(d)$. He used this connection to 
reduce the description of
 the possible singular spectra of products of 
 random matrices from  $SL(d,\comp)$ with given singular spectra
to the spectral problem for sums of Hermitian matrices. 
 Basically, Klyachko's observation is 
a connection between the convolution algebras of the Gelfand pairs
$(SL(d,\comp), SU(d))$ and $(SU(d)\ltimes H_d^0, SU(d))$.  
It is  closely related
to a similar correspondence between the convolutions of 
conjugation-invariant measures on a compact Lie group (here
  $SU(d)$) on one hand and  $Ad$-invariant measures on its Lie algebra (here $iH_d^0$) on the other, in terms of the so-called wrapping map 
(see \cite{DW} and  Remark \ref{wrapping}). Klyachko
 noticed in  \cite{Kl2}, but  did not not explain  that his
 connection   can 
be well expressed in terms of hypergroups; his proof goes via
random walks in the group $SU(d)$ and  relies on 
various identities between the spherical functions of 
$(SL(d,\comp), SU(d))$, the  characters of $SU(d)$ and  the Euclidean group
 $H_d^0$, as well as Poisson's summation formula.

The main purpose of the present paper  
is  to clarify and simplify Klyachko's approach \cite{Kl2}  by using so-called
 deformations of hypergroup convolutions by 
positive semicharacters as introduced in 
  \cite{V1},  \cite{V2}. 
Our description in particular
implies that the Banach algebra  $M_b(SL(d,\comp)\|SU(d))$
of  $SU(d)$-biinvariant signed measures on 
$SL(d,\comp)$
with total variation norm  is isometrically isomorphic to a certain
 Banach subalgebra of    $M_b(H_d^0)$, whereby probability measures are being
prserved.  We finally show  how this
leads to new proofs for  (known) limit theorems for $SU(d)$-biinvariant 
 random walks  on  $SL(d,\comp)$.
The explicit construction of this isomorphism runs as follows:
\begin{enumerate}
\item[\rm{(1)}] $SU(d)$ acts on  $H_d^0$ by conjugation
as a group of orthogonal transformations. The space $(H_d^0)^{SU(d)}$
of all orbits can be identified with 
 $$ C:=\{x\in\reel^d:
 \> x_1\ge x_2\ge\ldots\ge x_d,\> \sum_i x_i=0\},$$ 
where $x$ represents the ordered eigenvalues of a matrix in $H_d^0$.
The Banach algebra $M_b^{SU(d)}(H_d^0)$ of all $SU(d)$-invariant bounded
measures on $H_d^0$  can  be identified with the Banach algebra
$(M_b(C),*)$ of the associated orbit hypergroup  $C\simeq (H_d^0)^{SU(d)}$.
\item[\rm{(2)}] The polar decomposition of $SL(d, \b C)$ shows that 
the double coset space 
 $SL(d,\comp)//SU(d)$ may  also be identified with $C$ where now for 
$x\in C, \,\, e^x$
 corresponds to the singular spectrum $\sigma_{sing}(A) = \sigma(\sqrt{AA^*})$ of some
$A\in SL(d,\comp)$.
 This leads to a canonical Banach algebra  
isomorphism between $M_b(SL(d,\comp)\|SU(d))$ and the Banach algebra 
$(M_b(C),\bullet)$ of the  double coset hypergroup  
$C\simeq SL(d,\comp)//SU(d)$.
\item[\rm{(3)}] The characters of the hypergroups $(C,\bullet)$ and $(C,*)$
  are spherical functions of the corresponding symmetric spaces,
 and it is well-known 
that as functions on $C,$ they only 
differ by a known factor $J^{-1/2}(x)>0$ (see 
\cite{H1}). We prove
that such a connection between the characters of two  hypergroup
structures on $C$ implies that their convolutions are related
by a so-called hypergroup deformation. In particular, the supports
of the convolution products of two point measures are the same in both
hypergroups. This immediately implies the equivalence of 
the two spectral problems for Hermitian and unitary matrices, 
the main result
of  Klyachko \cite{Kl2}.
 The deformation isomorphism, however, is not isometric and not probability
preserving. To achieve this, a final correction is needed:
\item[\rm{(4)}] Instead of looking at  the Banach algebra
 $M_b^{SU(d)}(H_d^0)$ of all $SU(d)$-invariant measures on $H_d^0$,
we take a suitable exponential function $e_\rho:H_d^0\to]0,\infty[$ and 
observe that
the norm-closure of
$$\{e_\rho\cdot \mu:\> \mu\in M_b^{SU(d)}(H_d^0), \,\,
 \rm{supp}\,\mu \,\,\rm{compact}\}$$
is a Banach algebra which  turns out to be  isometrically
isomorphic to  $(M_b(C),\bullet)$.
\end{enumerate}
Putting all steps together  one obtains
the claimed probability preserving isometric isomorphism.
For $d=2$ this construction    reflects  the known close
 connection between the so-called 
Naimark hypergroup  $SL(2,\comp)//SU(2)$
 and  the  Bessel-Kingman hypergroup 
 $ (\reel^3)^{SO(3)}$;
 see \cite{BH},   \cite{V1}.

All results will be derived 
in  the general setting of a
 complex connected semisimple Lie group $G$ with finite center and 
 maximal compact subgroup $K$; in the special situation above, $(G,K)=(SL(d,\comp),SU(d))$.  For further results concerning the
associated hypergroup convolutions considered in this paper
we refer to \cite{DRW} and \cite{GS1}.

The paper is organized as follows: In Section 2 we collect  some relevant
 facts on commutative hypergroups  and their deformations.
In Section 3 we then use this deformation 
 to show how  the Banach algebras of all 
$K$-biinvariant bounded complex  measures on $G$ appear 
as   subalgebras of the  Banach algebra of 
 bounded complex measures on some Euclidean space. The final section is 
 devoted  to  probabilistic applications.

 \section{Deformations of commutative hypergroups}

In this section we give a quick introduction to hypergroups.
 We in particular prove a general result on deformations of
hypergroup convolutions in Proposition \ref{deform} which is crucial for step 
(3) of our construction. Moreover, step (4)  will 
be explained.

We first fix some  notations.
For a locally compact Hausdorff space  $X$, $M^+(X)$ denotes the 
 space of all positive 
Radon measures on $X$, and  $M_b(X)$ the Banach space of all bounded
regular complex Borel measures  with the total variation norm.
Moreover,  $M^1(X)\subset M_b(X) $ is the set of all probability measures,
$M_c(X)\subset M_b(X) $  the set of all measures with compact support, and 
 $\delta_x$  the point
measure in $x\in X$.
The spaces $C(X)\supset C_b(X)\supset C_0(X) \supset C_c(X)$ 
of continuous functions are defined as usual. For details on 
the following  we refer to \cite{BH} and  \cite{J}.

\begin{hypergroups} A hypergroup $(X,*)$ 
consists of a locally compact Hausdorff space $K$ and a convolution $*$
 on $M_b(X)$ 
such that $(M_b(X),*)$ becomes a Banach algebra, where  $*$ is
weakly continuous and probability preserving and  preserves compact
supports of measures.
 Moreover, there  exists an identity $e\in X$  
with  $\delta_e*\delta_x=
\delta_x*\delta_e=\delta_x$ for $x\in X$, as well as
  a continuous involution $x\mapsto\bar x$ on $X$ such that for $x, y\in X$,
$e\in supp(\delta_x*\delta_y)$ is equivalent to $ x=\bar y$,
and
$\delta_{\bar x} * \delta_{\bar y} =(\delta_y*\delta_x)^-$.
Here for $\mu\in M_b(X)$, the measure  $\mu^-$ is given by
$\mu^-(A)=\mu(A^-)$ for Borel sets $A\subset X$.

A hypergroup $(X,*)$ is called commutative if and only if so is the
convolution $*$. Hence, for a commutative hypergroup $(X,*)$ the triple   
$(M_b(X),*,.^-)$
is a commutative Banach-$*$-algebra with identity $\delta_e$. 
\end{hypergroups}

\begin{examples}\label{ex_1}
\begin{enumerate}
\item[\rm{(i)}] If $G$ is a  locally compact group, then $(G,*)$ is a
  hypergroup with the usual group convolution $*$.
\item[\rm{(ii)}]  Let $K$ be a compact subgroup of a locally compact group $G$.
Then 
$$M_b(G\|K):=\{\mu\in M_b(G):\> \delta_x*\mu*\delta_y=\mu\>\>\forall\>\> x,y\in K\}
$$
is a Banach-$*$-subalgebra of $M_b(G)$ with identity $dk$ where 
 $dk\in M^1(G)$ is the normalized Haar measure of $K$ embedded into $G$.
 Moreover, 
the double coset space $G/\!/K:=\{KxK:\> x\in G\}$ is  a locally compact
Hausdorff space, and  the canonical projection 
$p:G\to G/\!/K$ induces a  probability preserving, isometric isomorphism
$\,p: M_b(G\|K)\to M_b(G/\!/K)$ of Banach spaces by taking 
images of measures. The transport of
the convolution on  $M_b(G\|K)$ to $M_b(G/\!/K)$ via $p$ leads to a 
 hypergroup
structure $(G/\!/K, *)$ with identity $K\in G/\!/K$ and involution
$(KxK)^-:=Kx^{-1}K$, and  
$p$  becomes a  probability preserving,
 isometric isomorphism of Banach-$*$-algebras.
\item[\rm{(iii)}] Let $(V, \langle\,.\,,.\,\rangle)$ be a 
finite-dimensional Euclidean vector space and
 $K\subset O(V)$   a compact subgroup  of the orthogonal
group of $V,$ acting continuously on $V$. For $\mu\in M_b(V)$, 
denote the image measure of $\mu$ under $k\in K$ by $k(\mu).$
Then  the space of $K$-invariant measures
$$M_b^K(V):=\{\mu\in M_b(V):\> k(\mu)=\mu\>\>\forall \>\>
k\in K\}$$
is a Banach-$*$-subalgebra of 
 $M_b(V)$ (with the group convolution), with identity $\delta_0$.
Moreover,  the space $V^K:=\{K.x:\> x\in V\}$
 of all $K$-orbits in 
$V$ is   locally compact, and   the canonical projection 
$p:V\to V^K$  induces a probability preserving,
 isometric isomorphism
$p: M_b^K(V)\to M_b(V^K)$ 
of Banach spaces and an  associated
so-called orbit hypergroup  $(V^K, *)$ such that
$p$  becomes  a  probability preserving,
 isometric isomorphism of Banach-$*$-algebras. The involution 
on $(V^K, *)$ is given by $\overline{K.x} = -K.x$. 
\end{enumerate}
\end{examples}

We next collect some  data of a commutative hypergroup
 $(X,*)$. By a  result of
R. Spector, there exists a (up to normalization)  unique  
Haar measure $\omega\in M^+(X)$ which is characterized by
$\, \omega(f)=\omega(f_x)$ for all  $\,f\in C_c(X)$ and 
$x\in X,$ 
where we use the notation
\[f_x(y):=f(x*y):= \int_X f\> d(\delta_x*\delta_y).\]
Similar to the dual  of a locally compact abelian group, 
one defines the spaces
\begin{enumerate}\itemsep=-2pt
\item[\rm{(a)}] $\chi(X):=\{\alpha\in C(X):\> \alpha\ne 0,\>\>  \alpha(x*y)=
 \alpha(x)\alpha(y)\>\>\forall\>\> x,y\in X\}$;
\item[\rm{(b)}]  $X^*:=\{\alpha\in\chi(X):\> \alpha(\bar x)=
\overline{\alpha(x)}\>\>\forall\>\> x\in X\}$;
\item[\rm{(c)}]
$\widehat X:=X^*\cap C_b(X).$
\end{enumerate}
The elements of $X^*$ and  $\widehat X$ are called
semicharacters and characters, respectively.
All spaces above are locally compact Hausdorff spaces w.r.t. 
 the topology of 
  compact-uniform convergence.  

The Fourier transform on $ L^1(X, \omega)$ is defined by
\[\widehat f(\alpha):= \int_X f(x)\overline{\alpha(x)}\>
d\omega(x), \,\, \alpha\in \widehat X\,.\]
 Similar, the Fourier-Stieltjes transform of $\mu\in M_b(X)$ is
 defined by $\widehat \mu(\alpha) := \int_X \overline{\alpha(x)} d\mu(x),$ $\alpha\in \widehat X$.
 Both
 transforms
are injective, c.f. \cite{J}.
In the following, we 
 consider different hypergroup convolutions on $X$, and
 we write $\chi(X,*),\, \omega_*$ etc. in 
order to specify the relevant convolution.

\begin{examples}\label{ex_2}
\begin{enumerate}
\item[\rm{(i)}] Assume that in the situation of \ref{ex_1}(ii),
 $G/\!/K$ is commutative. Then a $K$-biinvariant function $\phi\in C(G)$ 
with $\phi(e)=1$ is by definition a spherical function of $(G,K)$ if 
$$\phi(g)\phi(h)=\int_K \phi(gkh)\> dk \quad\quad\forall \,\, g,h\in G.$$
Multiplicative functions $\alpha\in \chi(G/\!/K)$ are in one-to-one
correspondence with  spherical functions on $G$ via 
 $\alpha\mapsto \alpha\circ p$ for the  projection $p:G\to G/\!/K$.
In this way, the Fourier(-Stieltjes) transform on   $G/\!/K$
 corresponds to the  spherical Fourier(-Stieltjes) transform.
\item[\rm{(ii)}] In the situation of example \ref{ex_1}(iii),
the functions
 \[ \alpha_\lambda(K.x) =
 \int_K e^{i\langle \lambda,\,k.x\rangle} dk \quad (x\in V)\]
are  continuous multiplicative functions of the orbit hypergroup 
$(V^K,*)$ for  $\lambda\in V_\comp$, the complexification of $V$, where
$\alpha_\lambda \equiv \alpha_\mu$  holds if and only if
 $K.\lambda = K.\mu$. It is also well-known (see \cite{J}) that 
$\widehat{V^K}=\{\alpha_\lambda:\> \lambda\in V\}$. 
\end{enumerate}
\end{examples}

In  \cite{V1}, positive semicharacters were used to construct
deformed hypergroup convolutions. More precisely, the following 
 was proven there:

\begin{proposition} 
Let $\alpha_0\in X^*$
be a positive semicharacter on the commutative hypergroup $(X,*)$,
 i.e.,   $\alpha_0(x) >0$ for  $x\in X$.
Then the convolution
$$ \delta_x\bullet\delta_y :=
 \frac{1}{\alpha_0(x)\alpha_0(y)}\cdot \alpha_0(\delta_x * \delta_y)
\quad\quad\quad(x,y\in X)$$
 extends uniquely  to a bilinear, associative,
probability preserving, and
weakly  continuous convolution $\bullet$ on $M_b(X)$. Moreover, $(X,\bullet)$
becomes 
a commutative hypergroup with the identity and involution of $(X,*)$.
For  $\mu,\nu\in M_c(X)$, one has
\begin{equation}\label{deformconvo}
\alpha_0\mu\bullet\alpha_0\nu = \alpha_0(\mu*\nu).
\end{equation}
\end{proposition}

Note  that by Eq.(\ref{deformconvo}), the mapping $\mu\mapsto \alpha_0\mu$
establishes a canonical algebra isomorphism between $(M_c(X),*)$ and 
 $(M_c(X),\bullet)$ which usually -- when $\alpha_0$ is unbounded -- 
cannot be extended to $M_b(X)$.
The hypergroup  $(X,\bullet)$ is called the deformation of  $(X,*)$
w.r.t. $\alpha_0$. Clearly, many data of $(X,\bullet)$ can be expressed in
terms of $\alpha_0$ and corresponding  data of   $(X,*)$. 

\begin{proposition} In the above setting, we have\parskip=-1pt
\begin{enumerate}\itemsep=-2pt
\item[\rm{(i)}] $\omega_\bullet\,:= \alpha_0^2\omega_*$
is a Haar measure of $(X,\bullet)$.
\item[\rm{(ii)}] The mapping $M_{\alpha_0}: 
  \,\alpha\mapsto \alpha/\alpha_0\,$ is a homeomorphism (w.r.t. the
  compact-uniform topology) between $(X,*)^*$ and $ (X,\bullet)^*$, and also
  between $\chi(X,*)$ and $\chi(X,\bullet)$.  
\end{enumerate}
\end{proposition}

\begin{proof} For (i) and the first part of (ii) see \cite{V1}; 
the second part of (ii) is  analogous.
\end{proof}

We next turn to the following converse statement;  it will be
 crucial for  this paper:

\begin{proposition}\label{deform}
 Let $(X,*)$ and  $(X,\bullet)$ be commutative hypergroups
 on  $X$. Assume there is a positive semicharacter $\alpha_0$
 of $(X,*)$ such that the spaces of multiplicative continuous
functions for $(X,*)$ and
 $(X,\bullet)$ are related via
\[ \chi(X,\bullet) \,=\, \Big\{\frac{\alpha}{\alpha_0}\,: \, \alpha\in
 \chi(X,*)\Big\}.\]
Then $(X,\bullet)$ is 
the deformation of $(X,*)$ w.r.t. $\alpha_0$.
\end{proposition}

\begin{proof} 
Let $(X,\circ)$ denote the deformation of $(X,*)$ via $\alpha_0$. 
Take $\beta\in (X,\bullet)^\wedge$. Then by our assumption and
the proposition above, $\beta$ is multiplicative w.r.t. $\circ$ as well,
and the Fourier-Stieltjes 
transforms
 of $\delta_x\circ\delta_y$ and $\delta_x\bullet\delta_y$
 w.r.t. $(X,\bullet)$ satisfy 
\begin{eqnarray*} (\delta_x\circ\delta_y)^\wedge(\beta) \,&=&\, \int_X \overline{\beta(z)}
  d(\delta_x\circ\delta_y)(z)\,=\,\overline{\beta(x\circ y)} \\
  &=&\,
  \overline{\beta(x)} \cdot\overline{\beta(y)} \,=\,
  (\delta_x\bullet\delta_y)^\wedge(\beta).
\end{eqnarray*}
By the injectivity of the Fourier-Stieltjes-transform on $M_b(X)$, we obtain
$\delta_x\circ\delta_y = \delta_x\bullet\delta_y$. Thus the convolutions
of $(X,\bullet)$ and $(X,\circ)$ coincide, and so do the involutions,
because they are uniquely determined by the convolutions.
\end{proof}

\begin{example}
It is well-known (see \cite{BH},\cite{H1},\cite{J},\cite{Ko}) that
the double coset hypergroup $SL(2,\comp)//SU(2)$ may be realized as 
  hypergroup  $(X=[0,\infty[,\bullet)$  with
 multiplicative functions
$$\beta_\lambda(x):= \frac{\sin \lambda x}{\lambda\cdot \sinh x}
\quad\quad\quad (\lambda\in \comp).$$
Via the correspondence of $\beta_\lambda$ with $\lambda$, we have
   $\chi(X,\bullet) \simeq\comp$.
On the other hand, the orbit hypergroup  $(\reel^3)^{SO(3)}$ 
 may be realized as 
 the  Bessel-Kingman 
hypergroup  $(X=[0,\infty[,*)$ with
  multiplicative functions 
$$\alpha_\lambda(x)= \frac{\sin \lambda x}{\lambda x}
\quad\quad\quad (\lambda\in \comp);$$
see \cite{BH}, \cite{J}.
Clearly, $\alpha_i(x)= (\sinh x)/x\,$ is a positive semicharacter
on  $(X,*)$ abd  $\beta_\lambda = \alpha_\lambda/\alpha_i$ for
 $\lambda\in \comp$. Proposition \ref{deform}  implies
the known fact that  $(X,\bullet)$ is a
deformation of $(X,*)$, c.f. \cite{V1}.
\end{example}

We next consider further examples  which  explain step (4) 
 in the introduction; they
 are similar to a construction  in \cite{V2}. Let 
$(V, \langle\,.\,,\,.\,\rangle)$ be a 
Euclidean vector space of finite dimension $n$, 
 $K\subset O(V)$ a  compact subgroup  of the orthogonal group of $V$, and
  $(V^K, *)$   the associated  orbit hypergroup. Fix  $\rho\in V$ with 
$-\rho\in K.\rho$, 
and consider the exponential 
\[e_\rho(x):=e^{\langle \rho,x\rangle}\]
on $V$.  Let further
\[M_c^{\rho,K}(V):=\{e_\rho\mu:\>
 \mu\in M_c(V) \>\>K-\rm{invariant}\}.\]
The multiplicativity of  $e_\rho$
yields that $\, e_\rho\mu * e_\rho\nu =  e_\rho(\mu*\nu)$, where $*$ denotes the  group convolution on $V$. Hence 
$M_c^{\rho,K}(V)$ is a subalgebra of 
 the Banach-$*$-algebra $M_b(V)$, 
and its norm-closure 
$$M_b^{\rho, K}(V):= \overline{M_c^{\rho,K}(V)}$$
is  a Banach subalgebra of  $M_b(V)$.
Notice that  for $\rho\ne0$ this algebra is not closed under
the  involution on  $M_b(V)$; for instance, the
$n$-dimensional normal distribution $N_{\rho, I}$ with density
$ (2\pi)^{-n/2} e^{-|x-\rho|^2/2} $ 
is  contained in  $M_b^{\rho, K}(V)$
while this is not the case for  $N_{\rho, I}^*=N_{-\rho, I}$.
Nevertheless, we   prove that 
  $M_b^{\rho, K}(V)$ is isometrically  isomorphic as a Banach algebra 
to the Banach algebra of measures  of a suitable deformation of the 
orbit hypergroup   $(V^K, *)$. More precisely:

\begin{proposition}\label{euclidean-deform} Let $\rho\in V$ with $-\rho\in K.\rho$ and define
\[\alpha_0(K.x):=\int_K e_\rho(k.x)\>
  dk \quad (x\in V).\]
 Then  the following hold: \parskip=-1pt
\begin{enumerate}\itemsep=-1pt
\item[\rm{(i)}] $\alpha_0$ is   a positive semicharacter on  $(V^K, *).$ 
\item[\rm{(ii)}]  If  $(V^K,\bullet)$ is the deformation of 
  $(V^K, *)$  w.r.t.  $\alpha_0$, then the canonical projection 
$p: V\to V^K$ induces (by taking image measures) a
probability preserving isometric isomorphism of Banach algebras  
from  $M_b^{\rho, K}(V)$ onto  $M_b(V^K, \bullet)$.
\end{enumerate}
\end{proposition}

\begin{proof}
\begin{enumerate}
\item[\rm{(i)}]  $\alpha_0$ is  obviously  continuous and
  positive. For $x,y\in V$, we further have
\begin{eqnarray*}\alpha_0(K.x\,*\,K.y) &=& \int_{V^K} \alpha_0(K.z)\> 
d(\delta_{K.x}*\delta_{K.y})(K.z)
\\
&=& \int_K \int_K \int_K  e_\rho(k.(k_1.x\,+\,k_2.y)) \> 
 dk dk_1 dk_2\\
&=& \int_K \left(\int_K e_\rho(kk_1.x)\, dk_1\right)\cdot
 \left(\int_K e_\rho(kk_2.y)\,  dk_2\right)  dk\\
&=&\alpha_0(K.x)\,\alpha_0(K.y).
\end{eqnarray*}
Moreover, by the condition on $\rho$ above and the
 properties of the Haar measure  of $K$, 
$$\alpha_0(\overline{K.x}) = \alpha_0(-K.x)=\int_K e^{\langle -\rho,\,k.x\rangle}\>dk= \int_K e^{\langle \rho,\,k.x\rangle}\>
dk\,=\alpha_0(K.x).$$
Therefore,    $\alpha_0$ is a positive semicharacter on  $(V^K, *)$. 
\item[\rm{(ii)}] Consider the diagram
\begin{diagram}
M_c^{0,K}(V)& \rTo^{ \mu\>
\mapsto\> e_\rho\mu} &   M_c^{\rho,K}(V)\\
\dTo &&\dTo\\
(M_c(V^K), *)&
 \rTo^{ \nu\>
\mapsto\> \alpha_0\nu} &  (M_c(V^K), \bullet)
\end{diagram}
where  $M_c$ always stands for a space of  measures with
compact support, and  the vertical mappings  are obtained by
 taking image measures w.r.t.~$p$.
The mapping $p$  restricted to $M_c^{0,K}(V)$ and 
$M_c^{\rho,K}(V)$ is probability preserving and isometric. 
Moreover, for $\mu\in M_c^{0,K}(V)$ and $f\in C_c(V^K)$
 we obtain by
the $K$-invariance of $\mu$ and the definition of $\alpha_0$,
\begin{eqnarray*}
\int f\, dp(e_\rho\mu)&=& \int_{V} (f\circ p)\, e_\rho \> d\mu 
\\
&=&\int_K \int_{V} (f\circ p)(k.x) e_\rho(k.x)\> d\mu(x)\> dk
\\
 &=&\int_{V} (f\circ p)(x)\cdot\left(\int_K e_\rho(k.x)\>
   dk\right) d\mu(x)
\\
 &=& \int_{V} (f\circ p) (\alpha_0\circ p) \> d\mu
 = \int f\alpha_0\, dp(\mu),
 \end{eqnarray*}

which proves that the diagram  commutes. Therefore, as both 
horizontal mappings and the left vertical mapping are algebra isomorphisms,
 the  right vertical mapping is also a  
probability preserving isometric isomorphism of Banach algebras  
from  $M_c^{\rho, K}(V)$ onto  $M_c(V^K, \bullet)$.
A  continuity and density argument  shows that 
the Banach algebras  $M_b^{\rho, K}(V)$ and  
 $M_b(V^K, \bullet)$
are isomorphic via the probability preserving mapping $p$.
\end{enumerate}
\end{proof}

\section{ Biinvariant measures on complex noncompact semisimple Lie groups}

We  here identify
 the Banach algebra of all bounded 
 Borel measures on a 
connected  semisimple noncompact Lie group  with finite center, which are 
biinvariant under a  maximal compact subgroup, as a Banach  algebra of
 measures in some  Euclidean setting. 
For the general background we refer to \cite{H1}.

Let $G$ be a complex, noncompact connected 
semisimple  Lie group with finite center
and $K$  a maximal compact subgroup of $G$. 
Consider the corresponding  Cartan 
decomposition $\mathfrak g= \mathfrak k\oplus \mathfrak p$ of the Lie algebra of
$G$, and choose a  maximal abelian subalgebra $\mathfrak a\subseteq \mathfrak p$. $K$ acts on $\mathfrak p$ via the adjoint representation as
a group of orthogonal transformations with respect to the Killing form as
scalar product. Let further $W$ be the Weyl group of $K$, which acts on 
  $\mathfrak a$ as a finite reflection group, with  root system 
$\Sigma\subset \mathfrak a$. Here and lateron, $\mathfrak a$ is always 
identified with its dual $\mathfrak a^*$ via the Killing form, which we denote by  $\langle\,.\,,\,.\,\rangle$. We fix some Weyl
  chamber $\mathfrak a_+$ in $\mathfrak a$ and denote the associated system of
  positive roots by  $\Sigma^+$. The closed chamber
$C:=\overline{\mathfrak a_+}$ is a fundamental domain for the action 
of $W$ on $\mathfrak a$.
Later on we shall need the  half sum of roots,
$$\rho:= \sum_{\alpha\in\Sigma^+}\alpha\in \mathfrak a_+\,.$$

We now identify $C$ with the orbit hypergroup $(\mathfrak p^K, *)$ where each
$K$-orbit in $ \mathfrak p$ corresponds to its unique representative 
in $C\subset\mathfrak p $.  Then in view of example \ref{ex_2}(ii) above and  
Prop.~IV.4.8 of
\cite{H1}, the 
multiplicative continuous functions of $(C,*)$,
 considered as $K$-invariant functions on $\mathfrak p$, are given  by
\begin{equation}\label{har}
\psi_\lambda(x) =  \int_K e^{i\langle\lambda,\,k.x\rangle} dk 
\quad(x\in \mathfrak p)
\end{equation}
where $\lambda$ runs through the complexification $\mathfrak a_{\comp}$ 
 of $\mathfrak a$. 
Moreover, $\psi_\lambda\equiv \psi_\mu$ iff 
$\lambda$ and $\mu$ are in the same $W$-orbit. This is a special case of 
Harish-Chandra's integral formula for the spherical functions of a 
Cartan motion group. According to Theorem II.5.35 and Cor. II.5.36 
of \cite{H1},
they can also be written as
\begin{equation}\label{euclidean}
\psi_\lambda(x) = \frac{\pi(\rho)}{\pi(x)\pi(i\lambda)} 
\sum_{w\in W} (\det w) e^{i \langle \lambda,\,w.x\rangle}
\end{equation}
with the fundamental alternating polynomial
$$\pi(\lambda)=\prod_{\alpha\in \Sigma^+}\langle\alpha, \lambda\rangle.$$
On the other hand, $C$ can be identified
with  $G/\!/K$ where $x\in C$ corresponds
to the double coset $K(e^x)K$. According to this identification
 and the explicit
formula for the spherical
functions in  Theorem IV.5.7 of  \cite{H1},  the 
multiplicative continuous functions on the commutative double coset hypergroup
$(G/\!/K, \bullet)=(C,\bullet)$ are (as functions on $\mathfrak a$) given  by
\begin{equation}\label{hyperbolic}
\phi_\lambda(x) = \frac{\pi(\rho)\sum_{w\in W} (\det w) e^{\,i \langle \lambda, \, wx\rangle}}{\pi(i\lambda)\sum_{w\in W} (\det w) e^{\langle\rho,\,wx\rangle}}
\quad(x\in \mathfrak a)
\end{equation}
with $\lambda\in \mathfrak a_\comp$. Thus in particular,
\begin{equation}\label{zshg}
\phi_\lambda(x)=\frac{\psi_\lambda(x)}{\psi_{-i\rho}(x)} \quad\quad
\forall\,\, x\in \mathfrak a, \> \lambda\in \mathfrak a_\comp.
\end{equation}
Notice that $\psi_{-i\rho}$ is a positive semicharacter of $(C,*)$. By
Weyl's formula (\cite{H1}, Prop.~I.5.15.), 
$$\psi_{-i\rho}(x) =\, \prod_{\alpha\in \Sigma^+} 
\frac{\sinh
  \langle\alpha,x\rangle}{\langle\alpha,x\rangle}.$$
 The square of $\psi_{-i\rho}$ (called $J(x)$ in \cite{H1}) determines the ratio of the volume elements in $\mathfrak p$ and $G/K$.
Proposition \ref{deform} and  Eq.~(\ref{zshg}) show
  that  $(G/\!/K, \bullet)=(C,\bullet)$
is the deformation of the orbit hypergroup
 $(\mathfrak p^K, *) = (C,*)$ via
$\psi_{-i\rho}$. Moreover, we have 
$\, \psi_{-i\rho}(x) = \int_K e_\rho(k.x)dk\,$ with the half sum  $\rho\in C$ 
and $\,e_\rho(x)= e^{\langle \rho,x \rangle}$. As the condition
$-\rho\in K.\rho$ is satisfied,  Proposition \ref{euclidean-deform}
 further implies that the Banach algebra of measures
  of  $(G/\!/K, \bullet)=(C,\bullet)$ 
can be identified with 
$M_b^{\rho, K}(\mathfrak p)$,  which is the closure of
$$\{e_\rho\mu:\>
 \mu\in M_c(\mathfrak p) \>\>K-\rm{invariant}\}.$$
The claimed  isomorphism between 
 $M_b^{\rho, K}(\mathfrak p)$ and
$M_b(\mathfrak p^K, \bullet)\!\simeq\! (M_b(C),\bullet)\\ 
\simeq 
 M_b(G||K)$ is now  given by taking
image  measures w.r.t. the canonical projection 
$\mathfrak p\mapsto \mathfrak p^K$.

\begin{remark}\label{wrapping}
The algebra isomorphism $\,M_c(\mathfrak p^K,*)\to M_c(G/\!/K,\bullet), \\ \mu\mapsto \alpha_0\mu$ with the semicharacter $\alpha_0 = \psi_{-i\rho}$ is closely
related to the so-called wrapping map  for the compact Lie group $K$, see \cite{DW}.
In fact, as $G$ is complex, we have $\mathfrak p = i\mathfrak k$ in the 
Cartan decomposition $\mathfrak g = \mathfrak k \oplus \mathfrak p$. Thus
 $K \cong (K\times K)/K$ (on which $K$ acts by conjugation) is the dual symmetric space of $G/K$. The wrapping map $\Phi: M_b(\mathfrak k) \to 
M_b(K)$ is defined by
\[ \Phi(\mu)(f):= \mu(j\widetilde f\,), \quad f\in C(K)\]
where $\widetilde f(x) = f(\exp x)$ and $j: \mathfrak k \to \b R$ is the $K$-invariant extension of 
\[ j(x) = \prod_{\alpha\in \Sigma^+} 
\frac{\sin
  \langle\alpha,x\rangle}{\langle\alpha,x\rangle} = \psi_{-i\rho}(ix), \quad 
x\in \mathfrak a.\]
Notice  that $\|\Phi(\mu)\| \leq \|\mu\|$. As shown  in \cite{DW}, $\Phi$ 
is an algebra homomorphism from $K$-invariant 
measures in $M_b(\mathfrak k)$
to conjugation-invariant measures in $M_b(K)$. The proof thereof is based on Weyl's integration formula and Kirillov's character formula for compact groups. 
In contrast to our situation, $\Phi$ does not preserve positivity and can 
therefore not be associated with a hypergroup deformation.
\end{remark}

\begin{example}
The group $K=SU(d)$ is a maximal compact subgroup of the connected semisimple
Lie group  $G=SL(d,\comp)$. In the   Cartan decomposition  
$\mathfrak g= \mathfrak k \oplus\mathfrak p$ we obtain   $\mathfrak p$ as
 the additive group $H_d^0$ of
all Hermitian $d\times d$-matrices with trace $0$, on which
$SU(d)$ acts  by conjugation. 
Moreover, $\mathfrak a$
consists of all real diagonal matrices with trace 0 and will be identified
with
\[\{x = (x_1,\ldots,x_d)\in\reel^d:\> \sum_i x_i=0\}\]
on which the Weyl group acts as the symmetric group $S_d$ as usual.
We thus  may  take
 $$C:=\{x = (x_1,\ldots,x_d)\in\reel^d:
 \> x_1\ge x_2\ge\ldots\ge x_d,\> \sum_i x_i=0\}.$$
This set parametrizes the possible spectra of matrices from  $H_d^0$.
 A  system of positive roots corresponding to $C$ is $\Sigma^+ = \{e_i-e_j: 1\leq i<j\leq d\}$ where $e_1, \ldots, e_d$ 
denotes the standard basis of $\reel^d$. The root system is of type $A_{d-1}$.
In order to describe the probability preserving isometric isomorphism
 for this example explicitly, we realize the canonical
projections $\,q:SL(d,\comp)\to SL(d,\comp)//SU(d)\simeq C$
and $\,p: H_d^0 \to (H_d^0)^{SU(d)}\simeq C$ as follows:
For $A\in H_d^0 $, define  $p(A):= \sigma(A) \in C$ as the tuple of
eigenvalues of $A$, ordered by size. 
For $B\in SL(d,\comp)$  define $q(B)$ as the element  $x\in C$ such
that  the singular spectrum $\sigma_{sing}(B)=\sigma(\sqrt{BB^*})$ of $B$
is  $e^x:=(e^{x_1}\ge \ldots \ge e^{x_d})$.
We have
$$\rho=\sum_{\alpha\in\Sigma^+}\alpha \,= \left(d-1,d-3,
d-5, \ldots,-d+3,-d+1\right)
\in \mathfrak a_+$$
which implies that the Banach algebra  $M_b^{\rho, K}(\mathfrak p)$
above may be described  as the closure of
$$\{e_\rho\mu:\>
 \mu\in M_c(H_d^0) \>\>SU(d)-\rm{invariant}\},$$
with  
$$e_\rho(A)= 
\exp\bigl(tr[A\cdot {\rm diag}(d-1,d-3, d-5, \ldots,-d+3,-d+1)]\bigr).$$
As discussed above, the mappings $p,q$ induce (by taking images of
measures)  probability preserving isometric isomorphisms from the Banach algebras $M_b^{\rho, K}(\mathfrak p)$ and $M_b( SL(d,\comp)||SU(d))$   onto
$(M_b(C), \bullet)$ for the double coset convolution $\bullet$ on 
$C\simeq  SL(d,\comp)//SU(d)$.

Let us finally come back to the two spectral problems studied by Klyachko \cite{Kl2}. The  hypergroup deformation between the two convolutions $*$ and $\bullet$ on $C$ gives the natural explanation for the close connection
of these problems:
First, for fixed $x_1, x_2\in C,$ the probability measure 
$\delta_{x_1} * \delta_{x_2}$ is the
 distribution of possible spectra of sums $A_1 + A_2\in H_d^0$ where the
$A_i$ run through all matrices from $ H_d^0$ with $\sigma(A_i) = x_i$ 
($i=1,2$).
On the other hand, $d(\delta_{x_1}\bullet\delta_{x_2})(y)$ describes the  distribution of possible singular spectra $e^y$ of products $B_1B_2\in SL(d,\b C)$ where the
$B_i$ run through all matrices from $SL(d,\b C)$ with given singular
spectra $e^{x_i}$. As 
$\rm{supp}(\delta_{x_1}\bullet\delta_{x_2}) = \rm{supp}(\delta_{x_1} * \delta_{x_2})$ we obtain:

\begin{corollary} (c.f. \cite{Kl2}, Theorem B) For elements $x_1, x_2, y\in C$ the following are equivalent:\parskip=-2pt
\begin{enumerate}\itemsep=-2pt
\item[(i)] There exist matrices $\,A_1,  A_2\in H_d^0$ with given spectra
 $\sigma(A_i)=x_i\,$ such that $y=\sigma(A_1 + A_2)$. 
\item[(ii)] There exist matrices $\, B_i\in SL(d,\b C)$ with given singular
spectra
 $\sigma_{sing}(B_i) = e^{x_i}$ such that  $\,e^y= \sigma_{sing}(B_1B_2)$.
\end{enumerate} 
\end{corollary}

\end{example}

We now briefly discuss the further classical series of 
 complex, noncompact  connected simple  Lie groups with finite center (c.f.
Appendix C in  \cite{Kn}.)

\begin{series}\hfill
\begin{enumerate}\itemsep=-1pt
\item[\rm{(1)}] {\bf The $B_n$-case.}\, For $n\ge 2$ consider 
$G=SO(2n+1,\comp)$ with  the maximal compact subgroup $K=SO(2n+1)$.
In this case $\mathfrak a$ may be identified with $\reel^n$ with standard
basis $e_1,\ldots,e_n$, and we may choose 
$$C=\{x\in\reel^n:\> x_1\ge x_2\ge\cdots\ge x_n\ge0\}$$
and
$\Sigma^+=\{e_i\pm e_j:\> 1\le i<j\le n\}\cup\{e_i:\> 1\le i\le n\}$.
The Weyl group  $W$ is  isomorphic 
 to the semidirect product  $S_n \ltimes\ganz_2^n$, and 
$\rho=\left(2n-1,2n-3,\ldots,1\right).$
\item[\rm{(2)}] {\bf The $C_n$-case.}\, For $n\ge 3$ consider
$G=Sp(n,\comp)$ with  the maximal compact subgroup $K=Sp(2n+1)$.
In this case, again  $\mathfrak a=\reel^n$ with 
 $C$ and  $W$ as in the $B_n$-case.
A positive root system is
$$\Sigma^+=\{e_i\pm e_j:\> 1\le i<j\le n\}\cup\{2e_i:\> 1\le i\le n\},$$
and  $\rho=(2n,2n-2,\ldots,2)$. Comparing this with the 
 $B_n$-case, we see from  Eq.~(\ref{euclidean})  that the spherical functions 
$\psi_\lambda$ of the orbit hypergroups $\mathfrak p^K$ 
are the same in both cases.
The preceding results 
on hypergroup deformations therefore imply 
 that the double coset hypergroups $\,Sp(n,\comp)//Sp(2n+1)$ and 
$\,SO(2n+1,\comp)//SO(2n+1)$ are  
deformations of each other w.r.t. certain positive semicharacters.
\item[\rm{(3)}] {\bf The $D_n$-case.}\, For $n\ge 4$ consider 
$G=SO(2n,\comp)$ with the maximal compact subgroup $K=SO(2n)$.
In this case  $\mathfrak a=\reel^n$ with
$C=\{x\in\reel^n:\> x_1\ge x_2\ge\cdots\ge x_{n-1}\ge |x_n|\}$
and
$\Sigma^+=\{e_i\pm e_j:\> 1\le i<j\le n\}$. Thus
$\rho=(2n-2,2n-4,\ldots,2,0).$
\end{enumerate}
\end{series}

 \section{Limit theorems for biinvariant random walks}

In this final section we use our previous results to 
translate limit theorems 
 for random walks on $\mathfrak p$ into 
corresponding results for $K$-biinvariant random walks on $G$. We first sketch
the method.

\begin{method}\label{method}
Let $(Z_n)_{n\ge 1}$ be a sequence on independent $G$-valued random variables
with  distributions $\mu_n\in M^1(G||K)$ ($n\ge1$). Then
$(S_n:=Z_1 Z_2\ldots Z_n)_{n\ge 0}$ with  $S_0:=e$ forms a
$K$-biinvariant random walk on $G$. Using the canonical projection
 $q:G\mapsto G/\!/K\simeq C\,$ and the  associated 
 isomorphism $\,q: M_b(G||K)\mapsto (M_b(C),\bullet)$ of Banach algebras,
 we  see
 that $(q(S_n)))_{n\ge 0}$ is a Markov process on $C$ with initial
 distribution $\delta_0$ and  transition probabilities
$$P(q(S_n)\in A\> |\> q(S_{n-1})=x) \> =\> (\delta_x\bullet q(\mu_n))(A)$$
for $n\ge 1$, $x\in C$, and  Borel sets $A\subset C$. This means that
 $(q(S_n)))_{n\ge  0}$ is a
 random walk on the hypergroup $(C,\bullet)$ with transition
probabilities $(q(\mu_n))_{n\ge 1}\subset M^1(C)$. Clearly,  $q(S_n)$ has
distribution $$ q(\mu_1)\bullet\ldots \bullet q(\mu_n)= 
 q(\mu_1*\ldots *\mu_n).$$

On the other hand, when considering the canonical projection
 $p:\mathfrak p\to \mathfrak p^K\simeq C$, we find unique probability
measures $(\nu_n)_{n\ge 1}\in M_b^{K,\rho}(\mathfrak p)$ with $
p(\nu_n)= q(\mu_n)$ for $n\in\nat$. If $(X_n)_{n\ge1}$ is a 
 sequence on independent $\mathfrak p$-valued random variables with
 distributions $(\nu_n)_{n\ge 1}$, we get a random walk $(T_n:= \sum_{k=0}^n
 X_k)_{n\ge0}$ on $\mathfrak p$ whose projection $(p(T_n))_{n\ge0}$ is 
a random walk on $(C,\bullet)$ with the same transition probabilities as 
$(q(S_n))_{n\ge0}$, i.e.,  
 $(p(T_n))_{n\ge0}$ and $(q(S_n))_{n\ge0}$
have the same finite dimensional distributions and therefore admit
  the same  limit theorems. Hence,
 limit theorems for the random walk  $(T_n)_{n\ge0}$ on  $\mathfrak p$
may be translated into  limit theorems for $(q(S_n))_{n\ge0}$.

Here is an example. Assume that the $\mu_n=\mu$ are independent of $n$, and 
that for 
the associated
$\nu$ on $\mathfrak p$ the first moment  vector
\[m:=\int_{\mathfrak p} x\> d\nu(x) \in \mathfrak p\]
 exists. Then by
Kolmogorov's strong law, $T_n/n\to m$ almost surely as $n\to\infty$. In other words, 
  $\,T_n=nm+o(n)$ a.s..
 As   $p:\mathfrak p\to \mathfrak p^K\simeq C$
is  homogeneous of degree 1 and contractive 
(see below), we obtain
$p(T_n)=np(m)+o(n)$  and thus $p(T_n)/n\to p(m)$ a.s.. 
Therefore, $q(S_n)/n\to p(m)$ a.s.
 for $n\to\infty$.

This strong law for $(q(S_n))_{n\ge0}$  has the  computational
drawback that    $p(m)$ is described in terms of  
$\nu\in M_b^{K,\rho}(\mathfrak p)$.                                       
To overcome this, we introduce a suitable modified moment 
 function $m_1:C\to C$  such that for all 
probability measures
$\nu\in M_b^{K,\rho}(\mathfrak p)$ having first  moments,
\begin{equation}\label{momentenforderung}
p(m) = p(\int_{\mathfrak p} x\> d\nu(x))= \int_C m_1\> dq(\mu).
 \end{equation}
We shall show   that $m_1$ is determined uniquely by
(\ref{momentenforderung}) and give an 
 explicit formula.  This  allows 
 to compute the limiting constant $p(m)$
directly  on $C$.
\end{method}

\noindent
Let us now go into details. As indicated above, 
we need  
the following  properties of 
$p$:

\begin{lemma} Assume $\mathfrak p$ and $\mathfrak a$ carry the Euclidean norm
  given by the Killing form.\parskip=-1pt
\begin{enumerate}\itemsep=-2pt
\item[\rm{(i)}] $p$ is homogeneous of degree 1, i.e.
$\,p(tx)=tp(x)$ for all $t>0$ and $x\in 
\mathfrak p$.
\item[\rm{(ii)}] 
         $dist\,(K.x, K.y)\,=\, \|p(x)-p(y)\|\,$ for all $x,y\in\mathfrak p$.
\end{enumerate}
\end{lemma}

\begin{proof}
Part (i) is obvious. For part (ii), see \cite{H1}, Prop. I.5.18.
\end{proof}

We next define $m_1$. Motivated by 
(\ref{momentenforderung}) for  $\nu\in
M_b^{\rho,K}(\mathfrak p)$ with $p(\nu) = \delta_x$ for $x\in C$, we put
\begin{equation}\label{momentendef}
m_1(x):=\frac{1}{\psi_{-i\rho}(x)} 
\int_K k.x\cdot e^{\langle k.x,\rho\rangle}\> dk \quad (x\in \mathfrak
p).
 \end{equation}
Notice  here that $\psi_{-i\rho}(x)=\int_K e^{\langle k.x,\rho\rangle}\> dk$.
The  function $m_1$ is
obviously $K$-invariant, continuous, and satisfies
$ \|m_1(x)\|\le\|x\|$ for $x\in\mathfrak p$
with respect to the norm induced by the Killing form.

\begin{proposition}\label{m1inC}
\begin{enumerate}\itemsep=-1pt
\item[\rm{(i)}] We have $ m_1(x)\in C$ for  all $x\in \mathfrak p$.
\item[\rm{(ii)}] For all $x\in C$,
$$
m_1(x)= 
\frac{\sum_{w\in W} (\det w) e^{\,\langle wx,\rho\rangle}wx}
{\sum_{w\in W} (\det g) e^{\,\langle wx,\rho\rangle}}
- \sum_{\alpha\in \Sigma^+} \frac{\alpha}{\langle\alpha,\rho\rangle}.
$$
\end{enumerate}
\end{proposition}

\begin{proof} We first check that $m_1(x)\in\mathfrak a$ for 
 $x\in \mathfrak p$.
The definition of $m_1$ and the Harish-Chandra formula 
(\ref{har}) show that for $\zeta\in \mathfrak p$,
\begin{equation}\label{ableitungm1}
\langle m_1(x),\zeta\rangle\,=\,
\frac{1}{i\psi_{-i\rho}(x)}\cdot\partial_\zeta\psi_\lambda(x)\vert_{\lambda=
  -i\rho}
\end{equation}
where $\partial_\zeta$ is the derivative in direction $\zeta$ w.r.t. 
 $\lambda$.  The open Weyl chamber $\mathfrak
a_+$ corresponding to $C$ is an  orthogonal transversal 
 manifold for the adjoint
action of $K$ on $\mathfrak p$ (\cite{H1}, Ch.II, 3.4.(vi)). Therefore  the orthogonal
complement $\mathfrak a^\perp$ of $\mathfrak a$ in $\mathfrak p$ coincides
with the tangent space of the orbit $K.\rho$ in $\rho$. As 
$\lambda\mapsto \psi_{-i\lambda}(x)$ is constant on $K$-orbits, 
$\,\partial_\zeta \psi_{\lambda}(x)\vert_{\lambda=-i\rho} =0$ for 
$\zeta\in \mathfrak a^\perp$. 
Hence $m_1(x) \in (\mathfrak a^\perp)^\perp = \mathfrak a$.
In order to check $m_1(x)\in C$, we recall from Ch. 3 of \cite{GB} that
\begin{eqnarray*}\label{charC}
 C\,&=&\,\{x\in \mathfrak a:\> d(x,\rho)\leq d(x,w\rho)\quad\forall \,w\in
 W\}\\
 &=& \,\{x\in \mathfrak a:\> \langle x,\rho\rangle\geq \langle x,w\rho\rangle
   \quad\forall \,w\in W\}.
\end{eqnarray*} 
On the other hand, an elementary  rearrangement inequality 
 (Theorem 368
of \cite{HL}) shows that for $z\in \mathfrak p$ and $w^\prime\in W$,
\[ \sum_{w\in W} \langle z,\, w\rho\rangle\> e^{\langle z,\, w\rho\rangle}
\,\geq\,\sum_{w\in W} \langle z,\, ww^\prime\rho\rangle\> e^{\langle z,
  w\rho\rangle}.\]
As $\psi_{-i\rho} >0$, we obtain  for $w^\prime\in W$ that
\begin{eqnarray*}
\langle m_1(x),\rho\rangle \,&=&\, \frac{1}{\psi_{-i\rho}(x)}
    \int_{K} \langle k.x,\rho\rangle\> e^{\langle
  k.x,\rho\rangle} dk \\
 &=&\, \frac{1}{|W|\,\psi_{-i\rho}(x)}\int_{K}\sum_{w\in W}  \langle k.x,\,w\rho\rangle\> e^{\langle
  k.x,\,w\rho\rangle} dk\\ 
 &\geq &\, \frac{1}{|W|\,\psi_{-i\rho}(x)}\int_{K}\sum_{w\in W}  \langle k.x,\,ww^\prime\rho\rangle\> e^{\langle
  k.x,w\rho\rangle} dk\\   
 &=&\, \langle m_1(x),w^\prime\rho\rangle.
\end{eqnarray*}
Eq.~(\ref{charC}) now shows that $m_1(x)\in C$. This completes the proof of
 Part (i). Part (ii)  follows 
 from Eqs.~(\ref{ableitungm1}) and 
  (\ref{euclidean}) for the spherical functions $\psi_\lambda$
 on $\mathfrak a$.
\end{proof}

\begin{moments}\label{moments}
Let $\nu\in M_b^{K,\rho}(\mathfrak p)$ be a probability measure and $r>0$.
We say that $\nu$ admits $r$-th moments if
 $\int_{\mathfrak p} \|x\|^r \> d\nu(x)<\infty$, or
equivalently, if  
 $\int_{\mathfrak p} |\langle \xi, x\rangle|^r \> d\nu(x)<\infty$
 for all  $\xi\in\mathfrak p$.
This  condition can be  translated into a corresponding
 condition for  $p(\nu)\in M^1(C)$. In fact, as $K$ acts
on $\mathfrak p$ as a group of orthogonal transformations,  
$$ \int_{\mathfrak p} \|x\|^r \> d\nu(x)=
 \int_{\mathfrak p} \|p(x)\|^r \> d\nu(x)= \int_C  \|y\|^r \> dp(\nu)(y)
\in [0,\infty].$$
Therefore, $\nu\in M_b^{K,\rho}(\mathfrak p)$
admits $r$-th moments if and only if $p(\nu)$ admits $r$-th
moments.        Proposition \ref{m1inC} and 
the estimate $\, \|m_1(x)\|\le\|x\|$ for $x\in\mathfrak p$
 show that                 that
this condition for  $r\ge1$ implies that the              modified 
moment vector  $\int_C  m_1(y) \> dp(\nu)(y)\in C$ exists. 
Moreover, as
$$\int_{\mathfrak p} x\> d\nu(x)= 
\int_C\int_K k.x \,\frac{ e^{\langle k.x,\rho\rangle}}{\psi_{-i\rho}(x)}
\> dk\> dp(\nu)(x),$$
we obtain from  Propos.~\ref{m1inC} that $\,\int_{\mathfrak p} x\> d\nu(x)\in C$
and hence, as claimed, 
$$p\Bigl(\int_{\mathfrak p} x\> d\nu(x)\Bigr)=\int_{\mathfrak p} x\> d\nu(x)
=\int_C m_1 \> dp(\nu)\,\in C. $$
\end{moments}

As an application, we derive a 
strong law of large numbers of Marcin\-kiewicz-Zygmund:

\begin{theorem}\label{strong-law}
Let $r\in]0,2[$ and  $\mu\in M^1(G||K)$ such that
 $q(\mu)\in M^1(C)$ admits $r$-th moments.  Let $(Z_n)_{n\ge 1}$ be 
a  sequence of i.i.d. $G$-valued $\mu$-distributed random variables.
Then 
\[\frac{1}{n^{1/r}}
\Bigl(q(Z_1\cdot Z_2\cdots Z_n)-n \cdot c  \Bigr) \longrightarrow 0 \quad a.s.\]
for $n\to\infty$, with
 $\,\displaystyle c= \int_C m_1  dq(\mu)\,$ in case $r\in [1,2[$, while $c\in C$ is arbitrary
for $r\in]0,1[$.
\end{theorem}

\begin{proof} Let $r\in [1,2[$.
By  Section \ref{moments},  
 the  $r$-th moment of the associated 
$\nu\in  M_b^{K,\rho}(\mathfrak p)$ with $p(\nu)=q(\mu)$ exists. Let 
$(X_n)_{n\ge1}$ be i.i.d.  $\mathfrak p$-valued random variables
with distribution $\nu$, and $(T_n=X_1+\ldots +X_n)_{n\ge 0}$ the associated random
walk on $\mathfrak p$.
The classical  Marcinkiewicz-Zygmund law 
 (Theorem 5.2.2 of \cite{CT}) yields that
 for all $\xi\in\mathfrak p$, 
\[n^{-1/r}\bigl(\langle \xi, T_n\rangle -n\int_{\mathfrak p}
\langle \xi, x\rangle\> d\nu(x)\bigr)\to 0\quad a.s. \]
as $ n\to\infty.$
But this means that $\,T_n-n \int_{\mathfrak p} x\> d\nu(x) = o(n^{1/r})$ a.s.
and hence, 
by  \ref{m1inC}(i) and (\ref{momentenforderung}), 
$$p(T_n)-n\int_C m_1 dq(\mu)=  p(T_n)-n p\bigl(\int_{\mathfrak p} x\>
d\nu(x)\bigr)
=  o\bigl(n^{1/r}\bigr)\quad\quad a.s..$$
 As
$(q(Z_1Z_2\ldots Z_n))_{n\ge 0}$ and $p(T_n)_{n\ge 0}$ have the same 
finite-dimensional distributions, the claim follows.
The case  $r\in]0,1[$ is similar.
\end{proof}

For
 $(G,K)=(SL(d,\comp), SU(d))$, the mapping 
$q$
is given by $q(A):=(\ln a_1(A),\ldots,\ln a_d(A))$ where 
 $a_1(A)\ge  a_2(A)\ge\ldots \ge a_d(A)>0$ 
are the eigenvalues of $\sqrt{AA^*}$. We therefore obtain

\begin{corollary}\label{strong-law-spezial}
Let $r\in[1,2[$ and $\mu\in M^1(SL(d,\comp)|| SU(d))$ such that its 
projection 
$q(\mu)\in M^1(C)$ admits $r$-th moments.  Then, for each sequence
 $(Z_n)_{n\ge 1}$  of i.i.d. $SL(d,\comp)$-valued
 and $\mu$-distributed random variables,
\begin{eqnarray*}\frac{1}{n^{1/r}}
\Bigl(\bigl(\ln a_1(Z_1\cdot Z_2\cdots Z_n),\ldots,
 \ln a_d(Z_1\cdot Z_2\cdots Z_n)\bigr))&\phantom{+}& \\
-\,n \cdot \int_C m_1 \> dq(\mu)  \Bigr)
 &\longrightarrow & \,0\quad a.s..
\end{eqnarray*} 
\end{corollary}

\begin{remark}
The strong laws \ref{strong-law} and 
 \ref{strong-law-spezial} are in principle well-known especially for 
 $r=1$; see for instance \cite{BL} and references cited therein.
Nevertheless our  approach may be  of some interest, because it uses 
the close connection between  (biinvariant) random walks on $G$
and those on  $\mathfrak p$ in a simple explicit way.
\end{remark}

\end{article}
\end{document}